\documentclass[12pt]{article}
\usepackage{graphicx}
\usepackage{epsf, verbatim}
\usepackage{amssymb,amsmath}
\usepackage{lscape}
\usepackage{color}

\title{\textbf{More on the restricted almost unbiased Liu-estimator in Logistic regression}}
\small{\author{\textbf{Nagarajah Varathan }$^{1,2}$ and \textbf{Pushpakanthie Wijekoon}$^{3}$\\
\small{$^{1}$Postgraduate Institute of Science, University of Peradeniya, Sri Lanka}\\
\small{$^{2}$Department of Mathematics and Statistics, University of Jaffna, Sri Lanka}\\
\small{$^{3}$Department of Statistics and Computer Science, University of Peradeniya, Sri Lanka}\\
\small{\textit{email:~varathan10@gmail.com,~pushpaw@pdn.ac.lk}}}

\date{}
\setcounter{secnumdepth}{5} \setlength{\textheight}{8.7in}
\setlength{\textwidth}{6.8in} \setlength{\topmargin}{-0.5in}
\setlength{\headheight}{0 in} \setlength{\headsep}{0.6in}
\setlength{\oddsidemargin}{-0.03in} \setlength{\parindent}{1pc}

\numberwithin{equation}{section}
\begin{document}
\maketitle \pagestyle{myheadings}
\begin{abstract}
To address the problem of multicollinearity in the logistic
regression model, in this paper we propose a new estimator called
Stochastic restricted almost unbiased logistic Liu-estimator
(SRAULLE) when the prior information is available in the form of
stochastic linear restrictions. A Monte Carlo simulation study was
carried out to compare the performance of the proposed estimator
with some existing estimators in the scalar mean squared error
(SMSE) sense.
Finally, a real data example was given to appraise the performance of the estimators.
\end{abstract}
\fontsize{12}{20}
\section{Introduction}

The maximum likelihood estimation technique is commonly used method
to estimate the parameters of the logistic regression model. The
multicollinearity severely affects the varinace of the estimates of
parameters in the logistic regression model. As a result model
produces inefficient estimates. To overcome this issue, several
alternative estimators have been proposed in the literature. These
estimators were introduced mainly based on two types. The first type
of estimators are based only on the sample information and the
second type of estimators are based on the sample and priori
available information which may be in the form of exact or
stochastic linear restrictions. Logistic Ridge Estimator (LRE) by
Schaefer et al. (1984), the Principal Component Logistic Estimator
(PCLE) by Aguilera et al. (2006), the Modified Logistic Ridge
Estimator (MLRE) by Nja et al. (2013), the Logistic Liu Estimator
(LLE) by Mansson et al. (2012), the Liu-Type Logistic Estimator
(LTLE) by Inan and Erdogan (2013), the Almost Unbiased Ridge
Logistic Estimator (AURLE) by Wu and Asar (2016), the Almost
Unbiased Liu Logistic Estimator (AULLE) by Xinfeng (2015), and the
Optimal Generalized Logistic Estimator (OGLE) by Varathan and
Wijekoon (2017) are some of the first type of estimators proposed in
the literature. Under the second type of estimators, the Restricted
Maximum Likelihood Estimator (RMLE) by Duffy and Santner (1989), the
Restricted Logistic Liu Estimator (RLLE) by S¸iray et al. (2015),
the Modified Restricted Liu Estimator by Wu (2015), the Restricted
Logistic Ridge Estimator (RLRE) by Asar et al. (2016a), the
Restricted Liu-Type Logistic Estimator (RLTLE) by Asar et al.
(2016b), and the Restricted Almost Unbiased Ridge logistic Estimator
(RAURLE) by Varathan and Wijekoon (2016a) were introduced to improve
the performance of the logistic model when the exact linear
restrictions are available in addition to sample model. When the
restrictions on the parameters are stochastic, the Stochastic
Restricted Maximum Likelihood Estimator (SRMLE) (Nagarajah and
Wijekoon, 2015), the Stochastic Restricted Ridge Maximum Likelihood
Estimator (SRRMLE) (Varathan and Wijekoon, 2016b), and the
Stochastic Restricted Liu Maximum Likelihood Estimator (SRLMLE)
(Varathan and Wijekoon, 2016c) were proposed in the literature. In
this research, following Xinfeng (2015), a new estimator namely,
Stochastic restricted almost unbiased logistic Liu Estimator
(SRAULLE) is proposed for the logistic regression model with the
presence of stochastic linear restrictions as prior information. The
organization of the paper is as follows. The model specification and
estimation are given in Section 2. Proposed estimators and their
asymptotic properties are discussed in Section 3. In Section 4, the
conditions for superiority of SRAULLE over some existing estimators
are derived with respect to mean square error (MSE) criterion. A
Monte Carlo simulation study is conducted to investigate the
performance of the proposed estimator in the scalar mean squared
error (SMSE) sense in Section 5. A numerical example is discussed in
Section 6. Finally, some conclusive remarks are given in Section 7.
\section{Model Specification and estimation}
Consider the logistic regression model
        \begin{eqnarray}
            y_{i} = \pi_{i} + \varepsilon_{i}, ~~i= 1,...,n
        \end{eqnarray}
        which follows Binary distribution with parameter $\pi_{i}$ as
            \begin{eqnarray}
            \pi_{i} = \frac{\exp(x_{i}'\beta)}{1+\exp(x_{i}'\beta)},
        \end{eqnarray}
        where $x_{i}$ is the $i^{th}$ row of $X$, which is an $n\times p$ data matrix with $p$ explanatory variables and $\beta$ is a $p\times 1$ vector of coefficients, $\varepsilon_{i}$ are independent with mean zero and variance $\pi_{i}(1-\pi_{i})$ of the response $y_{i}$. The Maximum likelihood estimate (MLE) of $\beta$ can be obtained as follows:
             \begin{eqnarray}
                \hat{\beta}_{MLE}= C^{-1}X'\hat{W}Z,
            \end{eqnarray}
            where $C = X'\hat{W}X$; $Z$ is the column vector with $i^{th}$ element equals $logit(\hat{\pi}_{i})+ \frac{y_{i}-\hat{\pi}_{i}}{\hat{\pi}_{i}(1-\hat{\pi}_{i})}$ and $\hat{W} = diag[\hat{\pi}_{i}(1-\hat{\pi}_{i})]$. Note that $\hat{\beta}_{MLE}$ is an unbiased estimate of $\beta$ and its covariance matrix is
            \begin{eqnarray}
            Cov (\hat{\beta}_{MLE}) = \{X' \hat{W} X \}^{-1}.
            \end{eqnarray}
The MSE and SMSE of $\hat{\beta}_{MLE}$ are
\begin{eqnarray}
            MSE[\hat{\beta}_{MLE}] &=& Cov[\hat{\beta}_{MLE}]+ B[\hat{\beta}_{MLE}]B'[\hat{\beta}_{MLE}]\\ \nonumber
                                   &=& \{X' \hat{W} X \}^{-1}\\ \nonumber
                                   &=& C^{-1} \nonumber
\end{eqnarray}
and
\begin{eqnarray}
            SMSE[\hat{\beta}_{MLE}] &=&  tr[C^{-1}]\\ \nonumber
\end{eqnarray}
When the multicollinearity presents in the logistic regression model
(2.1), many alternative estimators have been proposed in the
literature. Among those, in this research we consider the Logistic
Liu estimator (LLE) by Mansson et al., (2012) and the Almost
Unbiased Logistic Liu Estimator (AULLE) by Xinfeng (2015) under the
first type of estimators.

        \begin{eqnarray}
            LLE : \hat{\beta}_{LLE}&=& Z_{d}\hat{\beta}_{MLE};
            ~where~~
            Z_{d}= (C + I)^{-1}(C + dI), 0 < d < 1 \\
            AULLE : \hat{\beta}_{AULLE}&=& W_{d}\hat{\beta}_{MLE};
            ~where~~
            W_{d}= [I-(1-d)^{2}(C + I)^{-2}], 0 < d < 1
        \end{eqnarray}

The asymptotic properties of LLE:

        \begin{eqnarray}
            E[\hat{\beta}_{LLE}] &=& E[Z_{d}\hat{\beta}_{MLE}]\\ \nonumber
                                &=& Z_{d}\beta,
         \end{eqnarray}

        \begin{eqnarray}
                    D[\hat{\beta}_{LLE}]&=& Cov[\hat{\beta}_{LLE}]\\ \nonumber
                                         &=& Cov[Z_{d}\hat{\beta}_{MLE}]] \\ \nonumber
                                 &=& Z_{d}C^{-1}Z_{d}',
        \end{eqnarray}

Consequently, the bias vector and the mean square error matrix of
LLE are obtained as
        \begin{eqnarray}
            B[\hat{\beta}_{LLE}] &=& E[\hat{\beta}_{LLE}]- \beta\\ \nonumber
                                &=& [Z_{d} - I]\beta \\ \nonumber
                                &=& \delta_{1}, (say)
         \end{eqnarray}
and
        \begin{eqnarray}
            MSE[\hat{\beta}_{LLE}] &=& D[\hat{\beta}_{LLE}]+ B[\hat{\beta}_{LLE}]B'[\hat{\beta}_{LLE}]\\ \nonumber
                                &=& Z_{d}C^{-1}Z_{d}' + \delta_{1}\delta_{1}'\\ \nonumber
         \end{eqnarray}
         respectively.\\\\
The asymptotic properties of AULLE:

        \begin{eqnarray}
            E[\hat{\beta}_{AULLE}] &=& E[W_{d}\hat{\beta}_{MLE}]\\ \nonumber
                                &=& W_{d}\beta,
         \end{eqnarray}

        \begin{eqnarray}
                    D[\hat{\beta}_{AULLE}]&=& Cov[\hat{\beta}_{AULLE}]\\ \nonumber
                                         &=& Cov[W_{d}\hat{\beta}_{MLE}]] \\ \nonumber
                                 &=& W_{d}C^{-1}W_{d}',
        \end{eqnarray}

Then, the bias vector and the mean square error matrix of AULLE are
obtained as
        \begin{eqnarray}
            B[\hat{\beta}_{AULLE}] &=& E[\hat{\beta}_{AULLE}]- \beta\\ \nonumber
                                &=& [W_{d} - I]\beta \\ \nonumber
                                &=& \delta_{2},
         \end{eqnarray}
and
        \begin{eqnarray}
            MSE[\hat{\beta}_{AULLE}] &=& D[\hat{\beta}_{AULLE}]+ B[\hat{\beta}_{AULLE}]B'[\hat{\beta}_{AULLE}]\\ \nonumber
                                &=& W_{d}C^{-1}W_{d}' + \delta_{2}\delta_{2}'\\ \nonumber
         \end{eqnarray}
         respectively.
As an alternative technique to stabilize the variance of the estimator due to multicollinearity, one can use prior information, if available, in addition to the sample model (2.1) either as exact linear restrictions or stochastic linear restrictions. \\\\
Suppose that the following stochastic linear prior information is
given in addition to the general logistic regression model (2.1).
         \begin{eqnarray}
            h = H\beta + \upsilon; ~~ E(\upsilon)= \textbf{0},~~  Cov(\upsilon) = \Psi.
         \end{eqnarray}
         where $h$ is an $(q\times 1)$ stochastic known vector, $H$ is a $(q\times p)$ of full rank $q\leq p$ known elements and $\upsilon$ is an  $(q\times 1)$ random vector of disturbances with mean \textbf{0} and dispersion matrix $\Psi$, which is assumed to be known $(q\times q)$ positive definite matrix. Further, it is assumed that $\upsilon$ is stochastically independent of $\varepsilon$, i.e) $E(\varepsilon\upsilon')= 0$.\\\\

\vspace{3mm} In the presence of stochastic linear restrictions
(2.17) in addition to the logistic regression model (2.1), Nagarajah
and Wijekoon (2015) introduced the Stochastic Restricted Maximum
Likelihood Estimator (SRMLE).
 \begin{eqnarray}
            \hat{\beta}_{SRMLE} = \hat{\beta}_{MLE} + C^{-1}H'(\Psi + HC^{-1}H')^{-1} (h - H\hat{\beta}_{MLE})
  \end{eqnarray}
  The asymptotic properties of SRMLE:
  \begin{eqnarray}
  E(\hat{\beta}_{SRMLE}) &=& \beta,
  \end{eqnarray}

    \begin{eqnarray}
  Cov(\hat{\beta}_{SRMLE}) &=& C^{-1} - C^{-1}H'(\Psi + HC^{-1}H')^{-1}HC^{-1}\\ \nonumber
                           &=& (C + H'\Psi^{-1}H)^{-1}, \\ \nonumber
                           &=& R (say)\nonumber
  \end{eqnarray}
 and \begin{eqnarray}
            Bias[\hat{\beta}_{SRMLE}] = E[\hat{\beta}_{SRMLE}] - \beta= 0.\\ \nonumber
         \end{eqnarray}
  The MSE of SRMLE is
 \begin{eqnarray}
            MSE[\hat{\beta}_{SRMLE}] &=&  Cov(\hat{\beta}_{SRMLE}) + B[\hat{\beta}_{SRMLE}]B'[\hat{\beta}_{SRMLE}]\\ \nonumber
                                     &=& (C + H'\Psi^{-1}H)^{-1}\\ \nonumber
                                     &=& R
\end{eqnarray}
\\
%
%
%
\section{The Proposed  Estimator}
In this section, by replacing $\hat{\beta}_{MLE}$ by $\hat{\beta}_{SRMLE}$ in (2.8), we propose a new estimator which is called as the Stochastic restricted almost unbiased logistic Liu Estimator (SRAULLE) and defined as\\
\begin{eqnarray}
            \hat{\beta}_{SRAULLE} &=& W_{d}\hat{\beta}_{SRMLE}
                                   \end{eqnarray}
  where $W_{d}= [I-(1-d)^{2}(C + I)^{-2}], 0 < d < 1$. \\\\
The asymptotic properties of $\hat{\beta}_{SRAULLE}$ are
        \begin{eqnarray}
            E[\hat{\beta}_{SRAULLE}] &= & E[W_{d}\hat{\beta}_{SRMLE}] \\ \nonumber
                                    &=& W_{d}\beta,\\ \nonumber
         \end{eqnarray}

   \begin{eqnarray}
    D(Cov(\hat{\beta}_{SRAULLE})&=& Cov(\hat{\beta}_{SRAULLE})\\ \nonumber
                              &=& Cov(W_{d}\hat{\beta}_{SRMLE})\\ \nonumber
                              &=& W_{d}Cov(\hat{\beta}_{SRMLE})W_{d}'\\ \nonumber
                              &=& W_{d}R W_{d}', \nonumber
    \end{eqnarray}
 and
 \begin{eqnarray}
 Bias(\hat{\beta}_{SRAULLE}) &=& E[\hat{\beta}_{SRAULLE}] - \beta \\ \nonumber
                            &=&  [W_{d}-I]\beta \\ \nonumber
                            &=& \delta_{2}.
                            \end{eqnarray}
Consequently, the mean square error can be obtained as,
       \begin{eqnarray}
        MSE(\hat{\beta}_{SRAULLE}) &=& D(\hat{\beta}_{SRAULLE}) + Bias(\hat{\beta}_{SRAULLE})Bias(\hat{\beta}_{SRAULLE})'\\ \nonumber
                                  &=& W_{d}R W_{d}' + \delta_{2}\delta_{2}' \\ \nonumber
       \end{eqnarray}

 \section{Mean square error comparisons}
 When different estimators are available for the same parameter
 vector $\beta$ in the regression model one must solve the
 problem of their comparison. Usually, as a general measure, the
 mean square error matrix is used, and is defined by
   \begin{eqnarray}
   MSE(\hat{\beta},\beta) &=& E[(\hat{\beta}- \beta)(\hat{\beta}- \beta)'] \\ \nonumber
                          &=& D(\hat{\beta})+ B(\hat{\beta})B'(\hat{\beta})
   \end{eqnarray}
   where $D(\hat{\beta})$ is the dispersion matrix, and $B(\hat{\beta})= E(\hat{\beta})- \beta$ denotes the bias vector. \\\\
   The Scalar Mean Squared Error (SMSE) of the estimator $\hat{\beta}$ can be defined as
   \begin{eqnarray}
    SMSE(\hat{\beta},\beta) = trace[MSE(\hat{\beta},\beta)]
  \end{eqnarray}
  For two given estimators $\hat{\beta}_{1}$ and $\hat{\beta}_{2}$, the estimator $\hat{\beta}_{2}$ is said to be superior to $\hat{\beta}_{1}$ under the MSE criterion if and only if
  \begin{eqnarray}
    M(\hat{\beta}_{1},\hat{\beta}_{2}) = MSE(\hat{\beta}_{1},\beta) - MSE(\hat{\beta}_{2},\beta) \geq 0.
  \end{eqnarray}
\subsection{Comparison of SRAULLE with MLE}
 To compare the estimators $\hat{\beta}_{MLE}$ and
$\hat{\beta}_{SRAULLE}$, we consider their MSE differences as below:
   \begin{eqnarray}
        MSE(\hat{\beta}_{MLE}) - MSE(\hat{\beta}_{SRAULLE}) &=& \{D(\hat{\beta}_{MLE})+ B(\hat{\beta}_{MLE})B'(\hat{\beta}_{MLE})\} \\ \nonumber
                                                          & &- \{ D(\hat{\beta}_{SRAULLE}) + B(\hat{\beta}_{SRAULLE})B'(\hat{\beta}_{SRAULLE})\}\\ \nonumber
                                                          &=& C^{-1} - \{W_{d}R W_{d}' + \delta_{2}\delta_{2}'\} \\ \nonumber
                                                          &=& U_{1} - V_{1}
   \end{eqnarray}
   where $U_{1} = C^{-1}$ and $V_{1} = W_{d}R W_{d}' + \delta_{2}\delta_{2}'$. One can obviously say that $W_{d}R W_{d}'$ and $U_{1}$
   are positive definite matrices and $\delta_{2}\delta_{2}'$ is non-negative definite matrix. Further by Lemma 1 (see Appendix A), it is clear that $V_{1}$
   is positive definite matrix. By Lemma 2 (see Appendix A), if $\lambda_{\rm{max}}(V_{1}U_{1}^{-1}) < 1$, then $U_{1} - V_{1}$ is a positive definite
   matrix, where $\lambda_{\rm{max}}(V_{1}U_{1}^{-1})$ is the largest eigen value of $V_{1}U_{1}^{-1}$. Based on the above arguments, it can be concluded
   that the estimator SRAULLE is superior to MLE if and only if $\lambda_{\rm{max}}(V_{1}U_{1}^{-1}) < 1$.

\subsection{Comparison of SRAULLE with LLE}
 Consider the MSE differences of $\hat{\beta}_{LLE}$ and
$\hat{\beta}_{SRAULLE}$
   \begin{eqnarray}
        MSE(\hat{\beta}_{LLE}) - MSE(\hat{\beta}_{SRAULLE}) &=& \{D(\hat{\beta}_{LLE})+ B(\hat{\beta}_{LLE})B'(\hat{\beta}_{LLE})\} \\ \nonumber
                                                          & &- \{D(\hat{\beta}_{SRAULLE})+ B(\hat{\beta}_{SRAULLE})B'(\hat{\beta}_{SRAULLE})\}\\ \nonumber
                                                          &=& \{Z_{d}C^{-1}Z_{d}' + \delta_{1}\delta_{1}'\} - \{W_{d}R W_{d}' + \delta_{2}\delta_{2}'\} \\ \nonumber
                                                          &=& U_{2} - V_{2}
   \end{eqnarray}
   where $U_{2} = Z_{d}C^{-1}Z_{d}' + \delta_{1}\delta_{1}'$ and $V_{2} = W_{d}R W_{d}' + \delta_{2}\delta_{2}'$. One can easily say that $W_{d}R W_{d}'$ and $Z_{d}C^{-1}Z_{d}'$
   are positive definite matrices and $\delta_{1}\delta_{1}'$ and $\delta_{2}\delta_{2}'$ are non-negative definite matrices. Further by Lemma 1, it is clear that
   $U_{2}$ and $V_{2}$ are positive definite matrices. By Lemma 2, if
$\lambda_{\rm{max}}(V_{2}U_{2}^{-1}) < 1$, then $U_{2} - V_{2}$ is a
positive definite
   matrix, where $\lambda_{\rm{max}}(V_{2}U_{2}^{-1})$ is the largest eigen value of $V_{2}U_{2}^{-1}$. Based on the above results, it can be
   said
   that the estimator SRAULLE is superior to LLE if and only if $\lambda_{\rm{max}}(V_{2}U_{2}^{-1}) < 1$.
\subsection{Comparison of SRAULLE with AULLE}
 Consider the MSE differences of $\hat{\beta}_{AULLE}$ and
$\hat{\beta}_{SRAULLE}$
   \begin{eqnarray}
        MSE(\hat{\beta}_{AULLE}) - MSE(\hat{\beta}_{SRAULLE}) &=& \{D(\hat{\beta}_{AULLE})+ B(\hat{\beta}_{AULLE})B'(\hat{\beta}_{AULLE})\} \\ \nonumber
                                                          & &- \{D(\hat{\beta}_{SRAULLE}) + B(\hat{\beta}_{SRAULLE})B'(\hat{\beta}_{SRAULLE})\}\\ \nonumber
                                                          &=& \{W_{d}C^{-1}W_{d}' + \delta_{2}\delta_{2}'\} - \{W_{d}R W_{d}' + \delta_{2}\delta_{2}'\} \\ \nonumber
                                                          &=& W_{d}(C^{-1}-R)W_{d}' \\ \nonumber
                                                          &=& C^{-1}H'(\Psi + HC^{-1}H')^{-1}HC^{-1} \\ \nonumber
                                                          & > &  0
   \end{eqnarray}
   Since the above mean square error difference is positive definite, it can be
   concluded that SRAULLE is always superior than AULLE.
\subsection{Comparison of SRAULLE with SRMLE}
 Consider the MSE differences of $\hat{\beta}_{SRMLE}$ and
$\hat{\beta}_{SRAULLE}$
   \begin{eqnarray}
        MSE(\hat{\beta}_{SRMLE}) - MSE(\hat{\beta}_{SRAULLE}) &=& \{D(\hat{\beta}_{SRMLE}) +B(\hat{\beta}_{SRMLE})B'(\hat{\beta}_{SRMLE})\} \\ \nonumber
                                                          & &- \{ D(\hat{\beta}_{SRAULLE})+ B(\hat{\beta}_{SRAULLE})B'(\hat{\beta}_{SRAULLE})\}\\ \nonumber
                                                          &=& R - \{W_{d}R W_{d}' + \delta_{2}\delta_{2}'\} \\ \nonumber
                                                          &=& U_{3} - V_{3}
   \end{eqnarray}
where $U_{3} = R$ and $V_{3} = W_{d}R W_{d}' +
\delta_{2}\delta_{2}'$. It can be easily seen that $W_{d}R W_{d}'$
and $R$
   are positive definite matrices and $\delta_{2}\delta_{2}'$ is non-negative definite matrix. Further by Lemma 1, it is clear that
   $V_{3}$ is positive definite matrix. By Lemma 2, if
$\lambda_{\rm{max}}(V_{3}U_{3}^{-1}) < 1$, then $U_{3} - V_{3}$ is a
positive definite
   matrix, where $\lambda_{\rm{max}}(V_{3}U_{3}^{-1})$ is the largest eigen value of $V_{3}U_{3}^{-1}$. Based on the above results, it can be
   said
   that the estimator SRAULLE is superior to SRMLE if and only if $\lambda_{\rm{max}}(V_{3}U_{3}^{-1}) < 1$.
\\\\
According to the results obtained from above mean square error
comparisons it can be concluded that the proposed estimator SRAULLE
is always superior than AULLE. However, under certain conditions
SRAULLE performs well over MLE, LLE, and SRMLE with respect to the
mean square error sense.

\section{ A Simulation study}
To examine the performance of the proposed estimator; SRAULLE with
the existing estimators: MLE, LLE, AULLE and SRMLE in this section,
we conduct the Monte Carlo simulation study. The simulations are
based on different levels of multicollinearity; $\rho= 0.7$, 0.8,
0.9 and 0.99 and different sample sizes; $n=25$, 50, 75 and 100. The
Scalar Mean Square Error (SMSE) is considered for the comparison.
Following McDonald and Galarneau (1975) and Kibria (2003), the
explanatory variables are generated as follows:
  \begin{eqnarray}
    x_{ij} = (1- \rho^{2})^{1/2} z_{ij} + \rho z_{i,p+1} , i = 1, 2, ...,n ,~~ j = 1, 2,...,p
  \end{eqnarray}
  where $z_{ij}$ are independent standard normal pseudo- random numbers and $\rho$ is specified so that the theoretical correlation between any two explanatory variables is given by $\rho^{2}$. Four explanatory variables are generated using (5.1).
   The dependent variable $y_{i}$ in (2.1) is obtained from the Bernoulli($\pi_{i}$) distribution where $\pi_{i} = \frac{\exp(x_{i}'\beta)}{1+\exp(x_{i}'\beta)}$. The parameter values of $\beta_{1}, \beta_{2},..., \beta_{p}$ are chosen so that $\sum_{j=1}^{p}\beta_{j}^{2}= 1$ and $\beta_{1} = \beta_{2} = ... = \beta_{p}$.
   Following Asar et al. (2016b), Wu and Asar (2015) and Mansson et al. (2012), the optimum value of the biasing parameter $d$ can be obtained by minimizing SMSE value with respect to $d$. However, for simplicity in this paper we consider some selected values of $d$ in the range $0 < d <1$.
   Moreover, we consider the following restrictions. \\
  \begin{eqnarray}
  H = \left(
        \begin{array}{cccc}
          1 & -1 & 0 & 1 \\
          1 & 1 & -1 & 0 \\
          0 & 0 & 1 & -1 \\
        \end{array}
      \right) ,~~
      h = \left(
            \begin{array}{c}
              1 \\
              -2 \\
              1 \\
            \end{array}
          \right) ~\rm{and}~~~ \Psi = \left(
                           \begin{array}{ccc}
                             1 & 0 & 0 \\
                             0 & 1 & 0 \\
                             0 & 0 & 1 \\
                           \end{array}
                         \right)
 \end{eqnarray}
The simulation is repeated 1000 times by generating new pseudo-
random numbers and the simulated SMSE values of the estimators are
obtained using the following equation.
 \begin{eqnarray}
    \hat{SMSE}(\hat{\beta}^{*}) &=& \frac{1}{1000}\sum_{r=1}^{1000}(\hat{\beta}_{r}- \beta)'(\hat{\beta}_{r}- \beta)
  \end{eqnarray}
where $\hat{\beta}_{r}$ is any estimator considered in the $r^{th}$
simulation. The simulation results are displayed in Tables 5.1 - 5.3
(Appendix). As we observed from the theoretical results, the
proposed estimator SRAULLE is superior to AULLE in the mean square
error sense with respect to all the sample sizes $n=25$, 50, 75 ,
and 100 and all the $\rho=0.7$, 0.8, 0.9, and 0.99. From the Tables
5.1- 5.3, it is further noted that if the multicollinearity is very
high (for example $\rho\geq 0.9$) the proposed estimator SRAULLE is
a very good alternative to MLE, LLE, AULLE and SRMLE regardless of
the values of $n$ and $d$. However, the performance of LLE is
considerably good for very small $d$ values and moderate $\rho$
values. Moreover, as we expected, MLE has the worst performance in
all of the cases (having the largest SMSE values).
\section{Numerical example}
In order to check the performance of the new estimator SRAULLE, in
this section, we used a real data set, which is taken from the
Statistics Sweden website (http://www.scb.se/). The data consists
the information about 100 municipalities of Sweden. The explanatory
variables considered in this study are Population ($x1$), Number
unemployed people ($x2$), Number of newly constructed buildings
($x3$), and Number of bankrupt firms ($x4$). The variable Net
population change ($y$) is considered as response variable, which is
defined as
 \begin{eqnarray*}
 y = \left\{\begin{array}{ccc}
       1 & ; & \rm{if ~there ~is ~an~ increase~ in ~the ~population}; \\
       0 & ; & o/w.
     \end{array}\right.
 \end{eqnarray*}
 The correlation matrix of the explanatory variables $x1$, $x2$, $x3$, and $x4$ is displayed in Table 6.1.
 It can be noticed from the Table 6.1 that, all the pair wise correlations are very high (greater than 0.95).
 Hence a clear high multicollinearity exists in this data set.
  Further, the condition number being a measure of multicollinearity is obtained as 188 showing that there
  exists severe multicollinearity with this data set. Moreover, we
  use the same restrictions as in (5.2) for the prior information.

   The SMSE values of MLE, LLE, AULLE, SRMLE, and SRAULLE for some selected values of
   biasing parameter $d$ in the range $0 <d < 1$ are given in the Table 6.2.

   It can be clearly noticed from the Table 6.2 that the proposed estimator SRAULLE outperforms the
   estimators MLE, LLE, AULLE, and SRMLE in the SMSE sense, with respect to all the selected values of biasing parameter $d$ in the range $0 < d < 1$ except $d=0.01$.
   Further, SRAULLE is having better performance compared to AULLE for all the values of $d$.

 \begin{center}
 \textbf{Table 6.1}: The correlation matrix of the explanatory variables\\
\begin{tabular}{ccccc}
  \hline
     & $x1$ & $x2$ & $x3$ & $x4$ \\ \hline
  $x1$ & 1.000 & 0.998 & 0.971 & 0.970 \\
  $x2$ & 0.998 & 1.000 & 0.960 & 0.958 \\
  $x3$ & 0.971 & 0.960 & 1.000 & 0.987 \\
  $x4$ & 0.970 & 0.958 & 0.987 & 1.000 \\
  \hline
\end{tabular}
\end{center}
\section{Concluding Remarks}
In this research, we proposed the Stochastic restricted almost
unbiased logistic Liu estimator (SRAULLE) for logistic regression
model in the presence of linear stochastic restriction when the
multicollinearity problem exists. The conditions for superiority of
the proposed estimator over some existing estimators were derived
with respect to MSE criterion. Further, a numerical example and a
Monte Carlo simulation study were done to illustrate the theoretical
findings. Results reveal that the proposed estimator is always
superior to AULLE in the mean square error sense and it can be a
better alternative to the other existing estimators under certain
conditions.
\\\\\\

\textbf{\large{References}}
\begin{enumerate}
        \item [1.] Aguilera, A. M., Escabias, M., Valderrama, M. J., (2006). Using principal components for estimating logistic regression with high-dimensional multicollinear data. \emph{Computational Statistics \& Data Analysis} 50: 1905-1924.
        \item [2.] Asar, Y., Arashi, M., Wu, J.,(2016a). Restricted ridge estimator in the logistic regression model. \emph{Commun. Statist. Simmu. Comp.}. Online. DOI: 10.1080/03610918.2016.1206932
        \item [3.] Asar, Y., Eri\c{s}o\v{g}lu, M., Arashi, M., (2016b). Developing a restricted two-parameter Liu-type estimator: A comparison of restricted estimators in the binary logistic regression model. \emph{Commun. Statist. Theor. Meth}. Online. DOI: 10.1080/03610926.2015.1137597
        \item [4.] Duffy, D. E., Santner, T. J., (1989). On the small sample prosperities of norm-restricted maximum likelihood estimators for logistic regression models. \emph{Commun. Statist. Theor. Meth}. 18: 959-980.
        \item [5.] Inan, D., Erdogan, B. E., (2013). Liu-Type logistic estimator. \emph{Communications in Statistics- Simulation and Computation}. 42: 1578-1586.
         \item [6.] Kibria, B. M. G., (2003). Performance of some new ridge regression estimators.\emph{Commun. Statist. Theor. Meth}. 32: 419-435.
        \item [7.] Mansson, G., Kibria, B. M. G., Shukur, G., (2012). On Liu estimators for the logit regression model. \emph{The Royal Institute of Techonology, Centre of Excellence for Science and Innovation Studies (CESIS)}, Sweden, Paper No. 259.
        \item [8.] McDonald, G. C., and Galarneau, D. I., (1975). A Monte Carlo evaluation of some ridge type estimators. \emph{Journal of the American Statistical Association} 70: 407-416.
        \item [9.] Nja, M. E., Ogoke, U. P., Nduka, E. C., (2013). The logistic regression model with a modified weight function. \emph{Journal of Statistical and Econometric Method} Vol.2, No. 4: 161-171.
        \item [10.] Rao, C. R., Toutenburg, H., Shalabh and Heumann, C., (2008). \emph{Linear Models and Generalizations}.Springer. Berlin.
        \item [11.] Rao, C. R.,and Toutenburg, H.,(1995). \emph{Linear Models :Least Squares and Alternatives, Second Edition}. Springer-Verlag New York, Inc.
        \item [12.] Schaefer, R. L., Roi, L. D., Wolfe, R. A., (1984). A ridge logistic estimator. \emph{Commun. Statist. Theor. Meth}. 13: 99-113.
        \item [13.] \c{S}iray, G. U., Toker, S., Ka\c{c}iranlar, S., (2015).On the restricted Liu estimator in logistic regression model. \emph{Communicationsin Statistics- Simulation and Computation} 44: 217-232.
         \item [14.] Nagarajah, V., Wijekoon, P., (2015). Stochastic Restricted Maximum Likelihood Estimator in Logistic Regression Model. \emph{Open Journal of Statistics}. 5, 837-851. DOI: 10.4236/ojs.2015.57082
         \item [15] Varathan, N., Wijekoon, P., (2016a). On the restricted almost unbiased ridge estimator in logistic regression. \emph{Open Journal of Statistics}. 6, 1076-1084. DOI: 10.4236/ojs.2016.66087
         \item [16.] Varathan, N., Wijekoon, P., (2016b). Ridge Estimator in Logistic Regression under stochastic linear restriction. \emph{British Journal of Mathematics \& Computer Science}. 15 (3), 1. DOI: 10.9734/BJMCS/2016/24585
         \item [17.] Varathan, N., Wijekoon, P., (2016c). Logistic Liu Estimator under stochastic linear restrictions. \emph{Statistical Papers}. Online. DOI: 10.1007/s00362-016-0856-6
         \item [18.] Varathan, N., Wijekoon, P., (2017). Optimal Generalized Logistic Estimator, \emph{Communications in Statistics-Theory and Methods} , DOI: 10.1080/03610926.2017.1307406
         \item [19.] Wu, J., (2015). Modified restricted Liu estimator in logistic regression model. \emph{Computational Statistics}. Online. DOI: 10.1007/s00180-015-0609-3
         \item [20.] Wu, J., Asar, Y., (2016). On almost unbiased ridge logistic estimator for the logistic regression model. \emph{Hacettepe Journal of Mathematics and Statistics}. 45(3), 989-998.\\ DOI: 10.15672/HJMS.20156911030
        \item [21.] Wu, J., Asar, Y., (2015). More on the restricted Liu Estimator in the logistic regression model. \emph{Communications in Statistics- Simulation and Computation}. Online. \\DOI: 10.1080/03610918.2015.1100735
        \item [22.] Xinfeng, C., (2015). On the almost unbiased ridge and Liu estimator in the logistic regression model.\emph{International Conference on Social Science, Education Management and Sports Education}. Atlantis Press: 1663-1665.

\end{enumerate}
\newpage
\begin{landscape}
\textbf{\large{Appendix}}\\\\
\textbf{Lemma 1:} Let $A$ : $n\times n$ and $B$ : $n\times n$ such that $A > 0$ and $B \geq 0$. Then $A + B > 0$. (Rao and Toutenburg, 1995)\\\\
\textbf{Lemma 2:} Let the two $n\times n$ matrices $M > 0$ ,$N \geq
0$, then $M > N$ if and only if $\lambda_{\rm{max}}(NM^{-1}) < 1$.
(Rao et al., 2008)\\\\\\
\textbf{Table 5.1}: The estimated MSE values for different $d$ when $n = 25$
\\\\
\footnotesize{\begin{tabular}{llccccccccccccc}
  \hline
     &                  & d = 0.01 & d = 0.1 & d = 0.2 & d = 0.3 & d = 0.4 & d = 0.5 & d = 0.6 & d = 0.7 & d = 0.8 & d = 0.9 & d = 0.99 \\\hline
$\rho = 0.70$ & MLE     & 1.4798 & 1.4798 & 1.4798 & 1.4798 & 1.4798 & 1.4798 & 1.4798 & 1.4798 & 1.4798 & 1.4798 & 1.4798   \\
              & LLE     & 0.8052 & 0.8600 & 0.9233 & 0.9893 & 1.0579 & 1.1291 & 1.2030 & 1.2794 & 1.3585 & 1.4402 & 1.5160  \\
              & AULLE   & 1.2383 & 1.2780 & 1.3186 & 1.3552 & 1.3875 & 1.4153 & 1.4383 & 1.4564 & 1.4694 & 1.4772 & 1.4798  \\
              & SRMLE   & 1.0138 & 1.0138 & 1.0138 & 1.0138 & 1.0138 & 1.0138 & 1.0138 & 1.0138 & 1.0138 & 1.0138 & 1.0138  \\
              & SRAULLE & 0.8536 & 0.8798 & 0.9067 & 0.9309 & 0.9524 & 0.9708 & 0.9861 & 0.9981 & 1.0068 & 1.0120 & 1.0137  \\
  \hline
$\rho = 0.80$ & MLE     & 1.9817 & 1.9817 & 1.9817 & 1.9817 & 1.9817 & 1.9817 & 1.9817 & 1.9817 & 1.9817 & 1.9817 & 1.9817   \\
              & LLE     & 0.8781 & 0.9602 & 1.0566 & 1.1584 & 1.2658 & 1.3786 & 1.4968 & 1.6205 & 1.7497 & 1.8843 & 2.0101  \\
              & AULLE   & 1.4912 & 1.5699 & 1.6513 & 1.7253 & 1.7912 & 1.8481 & 1.8955 & 1.9329 & 1.9599 & 1.9763 & 1.9817  \\
              & SRMLE   & 1.1793 & 1.1793 & 1.1793 & 1.1793 & 1.1793 & 1.1793 & 1.1793 & 1.1793 & 1.1793 & 1.1793 & 1.1793  \\
              & SRAULLE & 0.9039 & 0.9479 & 0.9935 & 1.0350 & 1.0720 & 1.1041 & 1.1308 & 1.1518 & 1.1670 & 1.1762 & 1.1793  \\
  \hline
$\rho = 0.90$ & MLE     & 3.5707 & 3.5707 & 3.5707 & 3.5707 & 3.5707 & 3.5707 & 3.5707 & 3.5707 & 3.5707 & 3.5707 & 3.5707   \\
              & LLE     & 0.9334 & 1.0954 & 1.2945 & 1.5137 & 1.7530 & 2.0124 & 2.2919 & 2.5915 & 2.9112 & 3.2510 & 3.5740  \\
              & AULLE   & 1.9075 & 2.1522 & 2.4151 & 2.6625 & 2.8885 & 3.0881 & 3.2573 & 3.3924 & 3.4908 & 3.5506 & 3.5705  \\
              & SRMLE   & 1.5271 & 1.5271 & 1.5271 & 1.5271 & 1.5271 & 1.5271 & 1.5271 & 1.5271 & 1.5271 & 1.5271 & 1.5271  \\
              & SRAULLE & 0.8744 & 0.9702 & 1.0733 & 1.1703 & 1.2590 & 1.3375 & 1.4039 & 1.4570 & 1.4957 & 1.5192 & 1.5270  \\
  \hline
$\rho = 0.99$ & MLE     & 33.1595 & 33.1595 & 33.1595 & 33.1595 & 33.1595 & 33.1595 & 33.1595 & 33.1595 & 33.1595 & 33.1595 & 33.1595   \\
              & LLE     & 0.4893 & 1.2132 & 2.5451 & 4.4324 & 6.8751 & 9.8731 & 13.4266 & 17.5353 & 22.1995 & 27.4190 & 32.5914  \\
              & AULLE   & 1.2984 & 3.5907 & 7.2401 & 11.5879 & 16.2114 & 20.7438 & 24.8751 & 28.3515 & 30.9757 & 32.6065 & 33.1540  \\
              & SRMLE   & 2.4804 & 2.4804 & 2.4804 & 2.4804 & 2.4804 & 2.4804 & 2.4804 & 2.4804 & 2.4804 & 2.4804 & 2.4804  \\
              & SRAULLE & 0.2878 & 0.4482 & 0.7008 & 1.0004 & 1.3183 & 1.6294 & 1.9128 & 2.1510 & 2.3308 & 2.4425 & 2.4800  \\
  \hline
\end{tabular}}
\newpage
\textbf{Table 5.2}: The estimated MSE values for different $d$ when $n = 50$ \\\\
\footnotesize{\begin{tabular}{llccccccccccccc}
  \hline
     &                  & d = 0.01 & d = 0.1 & d = 0.2 & d = 0.3 & d = 0.4 & d = 0.5 & d = 0.6 & d = 0.7 & d = 0.8 & d = 0.9 & d = 0.99 \\\hline
$\rho = 0.70$ & MLE     & 1.0662 & 1.0662 & 1.0662 & 1.0662 & 1.0662 & 1.0662 & 1.0662 & 1.0662 & 1.0662 & 1.0662 & 1.0662   \\
              & LLE     & 0.6249 & 0.6600 & 0.7004 & 0.7422 & 0.7854 & 0.8300 & 0.8761 & 0.9236 & 0.9726 & 1.0230 & 1.0696  \\
              & AULLE   & 0.9318 & 0.9543 & 0.9770 & 0.9975 & 1.0154 & 1.0307 & 1.0434 & 1.0533 & 1.0604 & 1.0647 & 1.0662  \\
              & SRMLE   & 0.7173 & 0.7173 & 0.7173 & 0.7173 & 0.7173 & 0.7173 & 0.7173 & 0.7173 & 0.7173 & 0.7173 & 0.7173  \\
              & SRAULLE & 0.6334 & 0.6474 & 0.6617 & 0.6744 & 0.6856 & 0.6952 & 0.7031 & 0.7092 & 0.7137 & 0.7164 & 0.7173  \\
  \hline
$\rho = 0.80$ & MLE     & 1.5025 & 1.5025 & 1.5025 & 1.5025 & 1.5025 & 1.5025 & 1.5025 & 1.5025 & 1.5025 & 1.5025 & 1.5025   \\
              & LLE     & 0.7156 & 0.7739 & 0.8419 & 0.9134 & 0.9883 & 1.0667 & 1.1485 & 1.2337 & 1.3224 & 1.4145 & 1.5003  \\
              & AULLE   & 1.1894 & 1.2403 & 1.2926 & 1.3400 & 1.3819 & 1.4181 & 1.4481 & 1.4717 & 1.4887 & 1.4990 & 1.5024  \\
              & SRMLE   & 0.8807 & 0.8807 & 0.8807 & 0.8807 & 0.8807 & 0.8807 & 0.8807 & 0.8807 & 0.8807 & 0.8807 & 0.8807 \\
              & SRAULLE & 0.7111 & 0.7388 & 0.7671 & 0.7928 & 0.8155 & 0.8351 & 0.8513 & 0.8641 & 0.8733 & 0.8789 & 0.8807  \\
  \hline
$\rho = 0.90$ & MLE     & 2.8448 & 2.8448 & 2.8448 &2.8448 & 2.8448 & 2.8448 & 2.8448 & 2.8448 & 2.8448 & 2.8448 & 2.8448  \\
              & LLE     & 0.8150 & 0.9431 & 1.0989 & 1.2687 & 1.4526 & 1.6505 & 1.8626 & 2.0887 & 2.3288 & 2.5830 & 2.8239  \\
              & AULLE   & 1.6588 & 1.8370 & 2.0267 & 2.2037 & 2.3644 & 2.5056 & 2.6248 & 2.7198 & 2.7888 & 2.8307 & 2.8446  \\
              & SRMLE   & 1.2238 & 1.2238 & 1.2238 & 1.2238 & 1.2238 & 1.2238 & 1.2238 & 1.2238 & 1.2238 & 1.2238 & 1.2238  \\
              & SRAULLE & 0.7500 & 0.8216 & 0.8976 & 0.9684 & 1.0325 & 1.0888 & 1.1363 & 1.1741 & 1.2015 & 1.2182 & 1.2237  \\
  \hline
$\rho = 0.99$ & MLE     & 26.9632 & 26.9632 & 26.9632 & 26.9632 & 26.9632 & 26.9632 & 26.9632 & 26.9632 & 26.9632 & 26.9632 & 26.9632   \\
              & LLE     & 0.4165 & 1.0607 & 2.1916 & 3.7594 & 5.7641 & 8.2058 & 11.0845 & 14.4001 & 18.1527 & 22.3422 & 26.4863  \\
              & AULLE   & 1.2822 & 3.2802& 6.3089 & 9.8379 & 13.5461 & 17.1550 & 20.4289 & 23.1753 & 25.2440 & 26.5281 & 26.9589  \\
              & SRMLE   & 2.2638 & 2.2638 & 2.2638 & 2.2638 & 2.2638 & 2.2638 & 2.2638 & 2.2638 & 2.2638 & 2.2638 & 2.2638  \\
              & SRAULLE & 0.1830 & 0.3484 & 0.5958 & 0.8823 & 1.1824 & 1.4737 & 1.7377 & 1.9589 & 2.1255 & 2.2288 & 2.2635  \\
  \hline
\end{tabular}}

\newpage
\textbf{Table 5.3}: The estimated MSE values for different $d$ when $n = 75$ \\\\
\footnotesize{\begin{tabular}{llccccccccccccc}
  \hline
     &                  & d = 0.01 & d = 0.1 & d = 0.2 & d = 0.3 & d = 0.4 & d = 0.5 & d = 0.6 & d = 0.7 & d = 0.8 & d = 0.9 & d = 0.99 \\\hline
$\rho = 0.70$ & MLE     & 0.4775 & 0.4775 & 0.4775 & 0.4775 & 0.4775 & 0.4775 & 0.4775 & 0.4775 & 0.4775 & 0.4775 & 0.4775   \\
              & LLE     & 0.3791 & 0.3887 & 0.3977 & 0.4076 & 0.4177 & 0.4280 & 0.4384 & 0.4489 & 0.4595 & 0.4703 & 0.4801  \\
              & AULLE   & 0.4647 & 0.4669 & 0.4691 & 0.4710 & 0.4728 & 0.4742 & 0.4754 & 0.4763 & 0.4770 & 0.4774 & 0.4775  \\
              & SRMLE   & 0.3940 & 0.3940 & 0.3940 & 0.3940 & 0.3940 & 0.3940 & 0.3940 & 0.3940 & 0.3940 & 0.3940 & 0.3940  \\
              & SRAULLE & 0.3837 & 0.3854 & 0.3872 & 0.3888 & 0.3902 & 0.3913 & 0.3923 & 0.3930 & 0.3936 & 0.3939 & 0.3940  \\
  \hline
$\rho = 0.80$ & MLE     & 0.6770 & 0.6770 & 0.6770 & 0.6770 & 0.6770 & 0.6770 & 0.6770 & 0.6770 & 0.6770 & 0.6770 & 0.6770  \\
              & LLE     & 0.4802 & 0.4968 & 0.5156 & 0.5349 & 0.5545 & 0.5745 & 0.5948 & 0.6156 & 0.6367 & 0.6582 & 0.6778  \\
              & AULLE   & 0.6404 & 0.6467 & 0.6530 & 0.6586 & 0.6634 & 0.6675 & 0.6709 & 0.6736 & 0.6755 & 0.6766 & 0.6770  \\
              & SRMLE   & 0.5120 & 0.5120 & 0.5120 & 0.5120 & 0.5120 & 0.5120 & 0.5120 & 0.5120 & 0.5120 & 0.5120 & 0.5120 \\
              & SRAULLE & 0.4850 & 0.4896 & 0.4943 & 0.4984 & 0.5020 & 0.5050 & 0.5075 & 0.5095 & 0.5109 & 0.5117 & 0.5120  \\
  \hline
$\rho = 0.90$ & MLE     & 1.3107 & 1.3107 & 1.3107 & 1.3107 & 1.3107 & 1.3107 &1.3107 & 1.3107 & 1.3107 & 1.3107 & 1.3107  \\
              & LLE     & 0.6774 & 0.7261 & 0.7824 & 0.8409 & 0.9016 & 0.9646 & 1.0299 &1.0974 & 1.1671 & 1.2391 & 1.3058  \\
              & AULLE   & 1.1005 & 1.1356 & 1.1712 & 1.2032 & 1.2312 & 1.2552 & 1.2750 & 1.2906 & 1.3017 & 1.3085 & 1.3107  \\
              & SRMLE   & 0.8007 & 0.8007 & 0.8007 & 0.8007 & 0.8007 & 0.8007 & 0.8007 & 0.8007 & 0.8007 & 0.8007 & 0.8007  \\
              & SRAULLE & 0.6742 & 0.6953 & 0.7167 & 0.7359 & 0.7528 & 0.7673 & 0.7792 & 0.7886 & 0.7953 & 0.7993 & 0.8007  \\
  \hline
$\rho = 0.99$ & MLE     & 13.1308 & 13.1308 & 13.1308 & 13.1308 & 13.1308 & 13.1308 & 13.1308 & 13.1308 & 13.1308 & 13.1308 & 13.1308  \\
              & LLE     & 0.5461 & 0.9882 & 1.6442 & 2.4735 & 3.4762 & 4.6523 & 6.0017 & 7.5245 & 9.2207 & 11.0902 & 12.9211  \\
              & AULLE   & 1.6681 & 2.8645 & 4.3929 & 6.0178 & 7.6355 & 9.1563 & 10.5043 & 11.6171 & 12.4466 & 12.9580 & 13.1291 \\
              & SRMLE   & 2.0259 & 2.0259 & 2.0259 & 2.0259 & 2.0259 & 2.0259 & 2.0259 & 2.0259 & 2.0259 & 2.0259 & 2.0259  \\
              & SRAULLE & 0.2834 & 0.4630 & 0.6941 & 0.9407 & 1.1869 & 1.4187 & 1.6245 & 1.7945 & 1.9212 & 1.9994 & 2.0256  \\
  \hline
\end{tabular}}
\newpage
\textbf{Table 5.4}: The estimated MSE values for different $d$ when $n = 100$ \\\\
\footnotesize{\begin{tabular}{llccccccccccccc}
  \hline
     &                  & d = 0.01 & d = 0.1 & d = 0.2 & d = 0.3 & d = 0.4 & d = 0.5 & d = 0.6 & d = 0.7 & d = 0.8 & d = 0.9 & d = 0.99 \\\hline
$\rho = 0.70$ & MLE     & 0.4172 & 0.4172 & 0.4172 & 0.4172 & 0.4172 & 0.4172 & 0.4172 & 0.4172 & 0.4172 & 0.4172 & 0.4172   \\
              & LLE     & 0.3380 & 0.3450 & 0.3528 & 0.3608 & 0.3689 & 0.3771 & 0.3853 & 0.3937 & 0.4021 & 0.4107 & 0.4185  \\
              & AULLE   & 0.4078 & 0.4094 & 0.4111 & 0.4125 & 0.4137 & 0.4148 & 0.4157 & 0.4163 & 0.4168 & 0.4171 & 0.4172  \\
              & SRMLE   & 0.3423 & 0.3423 & 0.3423 & 0.3423 & 0.3423 & 0.3423 & 0.3423 & 0.3423 & 0.3423 & 0.3423 & 0.3423  \\
              & SRAULLE & 0.3350 & 0.3363 & 0.3375 & 0.3387 & 0.3396 & 0.3404 & 0.3411 & 0.3416 & 0.3420 & 0.3422 & 0.3423  \\
  \hline
$\rho = 0.80$ & MLE     & 0.5932 & 0.5932 & 0.5932 & 0.5932 & 0.5932 & 0.5932 & 0.5932 & 0.5932 & 0.5932 & 0.5932 & 0.5932  \\
              & LLE     & 0.4341 & 0.4477 & 0.4629 & 0.4785 & 0.4943 & 0.5104 & 0.5268 & 0.5435 & 0.5605 & 0.5777 & 0.5935  \\
              & AULLE   & 0.5662 & 0.5709 & 0.5755 & 0.5796 & 0.5832 & 0.5863 & 0.5888 & 0.5907 & 0.5921 & 0.5930 & 0.5932  \\
              & SRMLE   & 0.4469 & 0.4469 & 0.4469 & 0.4469 & 0.4469 & 0.4469 & 0.4469 & 0.4469 & 0.4469 & 0.4469 & 0.4469 \\
              & SRAULLE & 0.4278 & 0.4311 & 0.4344 & 0.4373 & 0.4398 & 0.4420 & 0.4437 & 0.4451 & 0.4461 & 0.4467 & 0.4469  \\
  \hline
$\rho = 0.90$ & MLE     & 1.1311 & 1.1311 & 1.1311 & 1.1311 & 1.1311 & 1.1311 &1.1311 & 1.1311 & 1.1311 & 1.1311 & 1.1311 \\
              & LLE     & 0.6264 & 0.6659 & 0.7113 & 0.7584 & 0.8070 & 0.8573 & 0.9092 &0.9627 & 1.0179 & 1.0746 & 1.1271  \\
              & AULLE   & 0.9785 & 1.0041 & 1.0300 & 1.0532 & 1.0736 & 1.0910 & 1.1053 & 1.1166 & 1.1246 & 1.1295 & 1.1311  \\
              & SRMLE   & 0.6958 & 0.6958 & 0.6958 & 0.6958 & 0.6958 & 0.6958 & 0.6958 & 0.6958 & 0.6958 & 0.6958 & 0.6958  \\
              & SRAULLE & 0.6086 & 0.6232 & 0.6381 & 0.6513 & 0.6630 & 0.6729 & 0.6811 & 0.6875 & 0.6921 & 0.6949 & 0.6958  \\
  \hline
$\rho = 0.99$ & MLE     & 10.8045 & 10.8045 & 10.8045 & 10.8045 & 10.8045 & 10.8045 & 10.8045 & 10.8045 & 10.8045 & 10.8045 & 10.8045  \\
              & LLE     & 0.5945 & 0.9872 & 1.5492 & 2.2435 & 3.0700 & 4.0288 & 5.1198 & 6.3430 & 7.6985 & 9.1863 & 10.6384  \\
              & AULLE   & 1.7707 & 2.7752 & 4.0155 & 5.3061 & 6.5735 & 7.7539 & 8.7932 & 9.6474 & 10.2821 & 10.6727 & 10.8032 \\
              & SRMLE   & 1.8820 & 1.8820 & 1.8820 & 1.8820 & 1.8820 & 1.8820 & 1.8820 & 1.8820 & 1.8820 & 1.8820 & 1.8820  \\
              & SRAULLE & 0.3585 & 0.5315 & 0.7427 & 0.9609 & 1.1742 & 1.3722 & 1.5462 & 1.6889 & 1.7948 & 1.8600 & 1.8817  \\
  \hline
\end{tabular}}
\newpage
\textbf{Table 6.2}: The SMSE values of estimators for the Numerical example \\\\
\footnotesize{\begin{tabular}{lcccccc}
  \hline
         & d = 0.01 & d = 0.1 & d = 0.2 & d = 0.3 & d = 0.4 & d = 0.5  \\\hline
 MLE     & 0.0009457555 & 0.0009457555 & 0.0009457555 & 0.0009457555 & 0.0009457555 & 0.0009457555    \\
 LLE     & 0.0009441630 & 0.0009443098 & 0.0009444729 & 0.0009446361 & 0.0009447993 & 0.0009449624   \\
 AULLE   & 0.0009457541 & 0.0009457543 &  0.0009457546 & 0.0009457548 & 0.0009457550 & 0.0009457551   \\
 SRMLE   & 0.0009445487 & 0.0009445487 & 0.0009445487 & 0.0009445487 & 0.0009445487 & 0.0009445487   \\
 SRAULLE & 0.0009445472 & 0.0009445475 & 0.0009445477 & 0.0009445480 & 0.0009445481 & 0.0009445483   \\
  \hline
         & d = 0.6 & d = 0.7 & d = 0.8 & d = 0.9 & d = 0.99   \\\hline
 MLE     & 0.0009457555 & 0.0009457555 & 0.0009457555 & 0.0009457555 & 0.0009457555     \\
 LLE     & 0.0009451256 & 0.0009452888 & 0.0009454521 & 0.0009456153 & 0.0009457622  \\
 AULLE   & 0.0009457553 &  0.0009457554 & 0.0009457554 & 0.0009457555 & 0.0009457555    \\
 SRMLE   & 0.0009445487 & 0.0009445487 & 0.0009445487 & 0.0009445487 & 0.0009445487    \\
 SRAULLE & 0.0009445484 & 0.0009445485 & 0.0009445486 & 0.0009445487 & 0.0009445487    \\
  \hline
\end{tabular}}

\end{landscape}

\end{document}